# Hybrid-Dimensional Finite Volume Discretizations for Fractured Porous Media

Ivar Stefansson[1], Inga Berre[1,2], Eirik Keilegavlen[1]

[1]Department of Mathematics, University of Bergen

[2]Christian Michelsen Research

Corresponding author: Ivar.Stefansson@uib.no

**Abstract** Over the last decade, finite volume discretizations for flow in porous media have been extended to handle situations where fractures dominate the flow. These discretizations have successfully been combined with the discrete fracture-matrix models to yield mass conservative methods capable of explicitly incorporating the impact of fractures and their geometry. When combined with a hybrid-dimensional formulation, two central concerns are the restrictions arising from small cell sizes at fracture intersections and the coupling between fractures and matrix. Focusing on these aspects, we demonstrate how finite volume methods effectively can be extended to handle fractures, providing generalizations of previous work. We address the finite volume methods applying a general hierarchical formulation, facilitating implementation with extensive code reuse and providing a natural framework for coupling of different subdomains. Furthermore, we demonstrate how a Schur complement technique may be used to obtain a robust and versatile method for fracture intersection cell elimination. We investigate the accuracy of the proposed elimination method through a series of numerical simulations in 3D and 2D. The simulations, performed on fractured domains containing permeability heterogeneity and anisotropy, also demonstrate the flexibility of the hierarchical framework.

# 1 Introduction

When fractures are present in a porous media, flow may be totally dominated by the fractured structure and significant complexity is added to the problem of numerical modelling. Two main concerns are complex geometries, which challenge the grid generation, and high contrasts between permeabilities and length scales of matrix and fractures, which challenge the discretization schemes. The most common modelling approaches may be classified according to the relative importance of the fracture and matrix contribution and to which extent the fracture geometry is honoured. In the one extreme, where large scale fractures dominate completely, the rock matrix is irrelevant, and the model may be restricted to the fracture network. This gives the discrete fracture network models (DFN), e.g. (Long, Remer, Wilson, & Witherspoon, 1982; Endo & Long, 1984). At the other end of the spectrum, the domain may contain a high number of statistically homogeneous fine-scale fractures. In this case, it may be more reasonable to consider their effective impact on the flow instead of considering each fracture explicitly. This leads to equivalent continuum models with upscaled, possibly anisotropic, permeability fields, as described in (Sahimi, 2011).

Compromises between the two extremes also exist. Both the dual-porosity, the dual-permeability and more general multi-continuum models introduce computational domains of simplified geometry representing fractures and matrix connected through so-called transfer functions. In the dual-porosity models, all flow is assumed to take place through the fractures, and the matrix is only considered to have fluid storage capacity (Barenblatt, Zheltov, & Kochina, 1960; Warren & Root, 1963). In dual permeability and multi-continuum models, e.g. (Jarvis, Jansson, Dik, & Messing, 1991; Preuss & Narasimhan, 1985), flow is also allowed through disconnected blocks of the matrix, but a connected fracture network is assumed to carry the global flow. This allows for a scale separation where the effect of the smaller fractures is upscaled to the matrix blocks.

In this work, we will use the discrete fracture-matrix model (DFM) as described in (Dietrich, et al., 2005). The idea is to treat the rock surrounding explicitly represented fractures as a porous medium and possibly account for smaller, non-explicit fractures by an upscaled matrix permeability, as in the equivalent continuum models. The explicitly represented fractures are accounted for individually, as in the DFN methods. Reviews of fracture network models may be found in (Berkowitz, 2002; Dietrich, et al., 2005).

Motivated by the fractures being very thin compared to typical length scales in the reservoir, a common approach to the explicit fracture modelling introduced by Kiraly (1979) is to model them as lower-dimensional inclusions with an assigned hydraulic aperture, see also (Alboin, Jaffré, Roberts, Wang, & Serres, 2000; Flemisch, et al., 2018). Combined with a conforming discretization, in which fractures are located on boundaries between matrix cells, the result is a hybrid representation of the fractures. This approach has been used in combination with several different discretization techniques, see (Flemisch, et al., 2018) and references therein, and is the one we will adhere to.

Finite volume (FV) methods are among the most widely used discretizations for simulation of flow in porous media. The popularity is explained by the methods' local mass conservation as well as computational efficiency and flexibility. For the simplest version, the two-point approximations, the computational cost is low and the implementation reasonably straightforward, making it the principal workhorse for a range of porous media flow problems.

When it comes to handling of fractures in FV discretizations, the two main approaches are distinguished by whether the grids are restricted by the fracture geometry. Recently, there has been considerable interest in the embedded DFM models (EDFM) as part of a FV discretization, where the restriction of matrix grid alignment with the fractures is circumvented, see e.g. (Hajibeygi, Karvounis, & Jenny, 2011; Moinfar, Varavei, & Sepehrnoori, 2014). While beneficial for the gridding, the coupling between the discretization of fractures and matrix is more involved than in the classical FV DFM models using conforming grids that are considered in this paper. Such methods include the vertex-centred method by Reichenberger et al. (2006), the cell-centre based one introduced by Karimi-Fard et al. (2004) and extended by Sandve et al. (2012), and the related method of Ahmed et al. (2015) using transfer functions for the fracture-matrix coupling. Other approaches include the vertex approximate gradient and hybrid FV schemes (Brenner, Groza, Guichard, Lebeau, & Masson, 2016) and the two-phase method (Monteagudo & Firoozabadi, 2004), and refer to Geiger & Matthäi (2014) for a review of numerical methods based on the DFM model.

Two issues for FV methods introduced by explicit and conforming representation of fractures are the fracture-matrix coupling and the handling of intersecting fractures. We propose two general approaches aimed at enhancing computational efficiency and simplifying the modelling and implementation with respect to the two challenges, at minimal loss of solution quality.

The fracture-matrix coupling is approached through a hierarchical partitioning of the domain into a number of subdomains of different dimension corresponding to matrix, fractures, and fracture intersections (Boon, Nordbotten, & Yotov, 2017). Discretization is then performed for each subdomain independently. This leads to a simple implementation, and allows for the application of different discretization schemes for different domains, at the price of certain restrictions on the coupling between subdomains.

The challenge related to handling of intersecting fractures is due to the comparatively small volumes of cells at the one- and zero-dimensional fracture intersections, and was recognized by Slough et al. (1999) and Granet et al. (2001). These cells impair system matrix condition numbers and, in the case of transport or multi-phase simulations, time step restrictions. One widely used remedy for FV schemes is the elimination of these cells following the approach introduced by Karimi-Fard et al. (2004) based on the Star-Delta transformation. The approachhas also been extended by Sandve et al. (2012). While producing satisfactory results in many cases, the procedure does not handle the intersection of fractures of highly varying permeability and some aspects of multiphase flow, as shown in (Stefansson, 2016; Walton, Unger, Ioannidis, & Parker, 2017).

To eliminate intersection cells while accounting for assigned permeabilities in the crossing fractures, we suggest a new approach to the elimination of intersections that accounts for the permeability of the intersection cell. By discretizing the domain with the intersection cells and then performing a Schur complement reduction to remove the degrees of freedom corresponding to the intersection, we obtain a reduced system possessing information on the permeability of the eliminated cells. Moreover, we show how the original elimination based on the Star-Delta transformation can be interpreted as the limiting case of the new scheme as the intersection permeability goes to infinity. The new technique also provides a natural way of back-calculating the pressure values at the intersections, as may be desirable e.g. for coupled problems. Finally, the technique is not restricted to FV methods and may in principle also be applied to other discretizations of the problem.

For a discretization to be suitable for DFM modelling it should be applicable to unstructured grids corresponding to geometrically complex fracture networks and strong parameter heterogeneity. The ability to handle strong anisotropy in permeability is also necessary in many situations. For a hybrid-dimensional FV method on a conforming grid, the geometrical complexity is mostly dealt with by the constrained meshing of the domain. The handling of parameter heterogeneity and, in particular, the anisotropy, is related to the particular FV method and flux approximation applied. These aspects will also be investigated in the current work.

The transport of heat or chemical species is oftentimes of primary interest in porous media modelling. Because of this, it is sensible to evaluate flow methods not only in terms of pressure fields, but also indirectly through transport simulations on the flow fields they produce. By examining accumulations of tracer, we aim at revealing details and differences that may be critical to the concentration distribution, but almost indiscernible in pressure comparisons.

The implementation used in this work is available in the open-source simulation tool PorePy (Keilegavlen, Fumagalli, Berge, Stefansson, & Berre, 2017). In addition to open source code, we provide documented run scripts for the results shown in Section 3 at www.github.com/pmgbergen/porepy. The exception to the above is the implementation of the Star-Delta elimination procedure, which is used for comparison with the Schur complement procedure, and for which we rely on the Matlab Reservoir Simulation Toolbox (Lie, 2016).

## 2 Models

As noted in the introduction, the DFM approach to fracture modelling acknowledges that fractures will strongly influence the behaviour of processes in porous media by explicit modelling of selected fractures. Any remaining fractures are considered homogenized and are represented by the matrix through average quantities. This motivates a structured framework for the domain description, as well as governing equations for both matrix and fractures. After the two have been introduced, we describe the FV discretizations used and elaborate on the general methodology for subdomain coupling and the approach to handle fracture intersections.

### 2.1 Computational grid

We employ the hierarchical subdomain approach to discrete fracture modelling and associated notation introduced by Boon, Nordbotten and Yotov (2017). The overarching idea is to use the constraints imposed by the fracture geometry to divide the $N$-dimensional domain $\Omega$ into subdomains $\Omega_i^d$ of dimension $d \leq N$. Specifically, for $N = 3$, we extract the matrix without fractures as a 3D subdomain, the fractures as 2D subdomains, fracture plane intersections as 1D subdomains and the intersection points of such lines as 0D subdomains. A similar cascade is defined for $N = 2$. As an example, consider the domain depicted in Figure 1, which consists of one subdomain of dimension three, two fracture subdomains of dimension two and one intersection subdomain of dimension one. 0D intersections were not included to ensure visual simplicity.

Each subdomain is meshed with the immersed lower-dimensional subdomains acting as constraints. We define a mesh size indicator $h_{min}$ to be the smallest cell diameter of the mesh of the highest dimension.

Although not required for the general hierarchical approach, we further assume the grids to match between subdomains of subsequent dimensions, so that each face in $\overline{\Omega}_i^d \cap \Gamma_{ij}$ corresponds to exactly one cell in $\Omega_j^{d-1}$. For most realistic fracture networks this demands unstructured grids, of which we limit ourselves to simplex grids in this work. For analytical purposes, however, it may also be useful to consider Cartesian grids. All descriptions in the following apply to simplex and Cartesian grids alike.

All subdomains are discretized separately before they are coupled two at a time by discretizations on the interfaces $\Gamma_{ij}$ between subdomains $i$ and $j$. This requires the coupling conditions to be local to the interfaces, but makes for a very flexible method in several ways. Firstly, we can use different discretizations in the different subdomains, as will be shown below. Secondly, this makes the conversion from a general porous medium discretization into a mixed-dimensional DFM method straightforward: An existing mono-dimensional method may be coupled by an appropriate coupling scheme. Thirdly, this also in principle facilitates the use of different physical models in different parts of the domain.

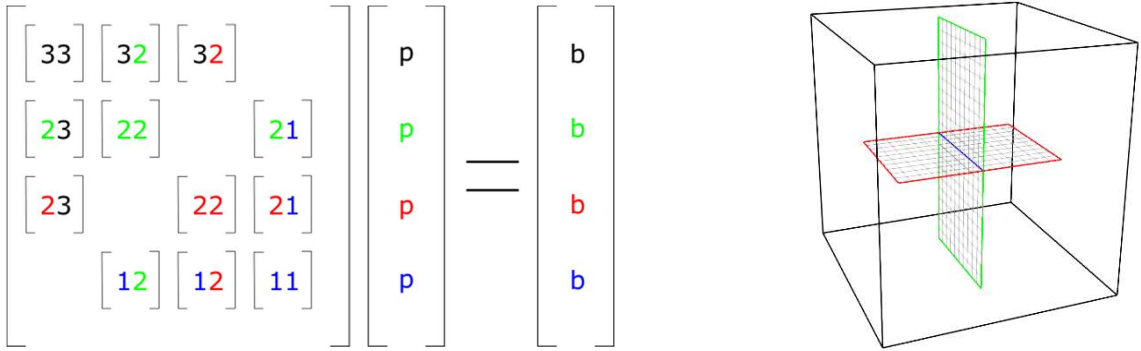

*Figure 1. To the left, a linear system matrix with off-diagonal subdomain coupling blocks and diagonal blocks for the internal subdomain discretizations. To the right, the domain with one, two and one subdomains of dimensions 3, 2 and 1, respectively. The correspondence between subdomains and matrix blocks is indicated through dimension numbers and colours.*

## 2.2 Governing equations

We solve the single-phase mixed-dimensional Darcy flow problem

$$\begin{aligned} \mathbf{u} + \mathbf{K}\nabla p &= 0, \\ \nabla \cdot \mathbf{u} - [\![\, \hat{\mathbf{n}} \cdot \hat{\mathbf{u}} \,]\!] &= f \end{aligned} \qquad \forall\, \Omega_i^d \subset \Omega. \qquad (1)$$

Here, $\mathbf{u}$ denotes the flux, $\mathbf{K}$ the permeability tensor, $p$ the pressure, $f$ the sources and sinks, and $\mathbf{n}$ the outward unit normal vector of a subdomain of dimension $d$. Quantities related to dimension $d+1$ are indicated using the hat notation $\hat{\cdot}$. Denoting the jump operator from the neighbouring higher-dimensional grids by $[\![\cdot]\!]$, the second term in the second equation of (1) enforces coupling of the subdomains. For the highest dimension $N$, we define $\hat{\mathbf{u}}$ to be zero.

We divide the boundary into parts where we prescribe Dirichlet and Neumann type boundary conditions, respectively:

$$\begin{aligned} p &= p_D & \text{on } \partial\Omega_{D,f}, \\ \mathbf{u} \cdot \mathbf{n} &= u_N & \text{on } \partial\Omega_{N,f}, \end{aligned} \qquad (2)$$

using subscript $f$ to indicate the flow problem.

We also model the advection of a passive tracer concentration $T$. For a fixed flux field and unitary porosity and density, the concentration is given according to

$$\frac{\partial T}{\partial t} + \nabla \cdot (T\mathbf{u}) - [\![\, \hat{T}\, \hat{\mathbf{n}} \cdot \hat{\mathbf{u}} \,]\!] = s, \qquad (3)$$

with $s$ denoting the tracer sources and sinks, and the corresponding boundary conditions

$$\begin{aligned} T &= T_D & \text{on } \partial\Omega_{D,t} \\ T\mathbf{u} \cdot \mathbf{n} &= J_N & \text{on } \partial\Omega_{N,t}. \end{aligned} \qquad (4)$$

## 2.3 Finite Volume Discretizations

We explore the subdomain coupling and intersection cell elimination in a FV framework. To obtain the formulation, we start out from the integral form of the mass balance of Equation (1) over a cell $c_i$. Using Gauss' divergence theorem on the left-hand side, we obtain

$$\int_{\delta c_i \backslash \Gamma} \mathbf{u} \cdot \mathbf{n} dS - \int_{\delta c_i \cap \Gamma} [\![ \hat{\mathbf{n}} \cdot \hat{\mathbf{u}} ]\!] dS = \int_{c_i} f dV, \qquad (5)$$

where $\mathbf{n}$ is the outward normal vector on the cell's boundary $\delta c_i$. Assuming the right-hand side source term to be known, and postponing the treatment of the coupling term to Section 2.4, we proceed to the approximation of the flux in the interior of the subdomain. Noting its dependency on the pressure as described by Equation (1), we aim at a relationship between the flux over a face $S_{ij}$ between cells $i$ and $j$ in terms of the pressures $p_k$ at the centres of the surrounding cells, i.e.,

$$\int_{S_{ij}} \mathbf{u} \cdot \mathbf{n} dS = -\int_{S_{ij}} \nabla p \cdot \mathbf{K} \cdot \mathbf{n} dS \approx \sum_k t_k p_k. \qquad (6)$$

The weights, or transmissibilities, $t_k$ incorporate geometry and permeability, and may be computed in different manners leading to different FV discretizations. In this work, we limit ourselves to two versions: the two-point and multi-point flux approximations described in the following subsections (TPFA and MPFA, respectively). Once computed, the transmissibilities are assembled into the discretization matrix $\mathbf{A}$, leading to a global system of equations of the form

$$\mathbf{Ap} = \mathbf{b}, \qquad (7)$$

where $\mathbf{b}$ incorporates boundary conditions and source and sink terms and $\mathbf{p}$ is the vector of cell centre pressures.

Note that even if we only solve for pressures, the fluxes may be readily back-calculated using the information in $\mathbf{A}$, which indeed is nothing else than a discretization and summation of the fluxes in the domain.

### 2.3.1 Two-Point Flux Approximation

To prevent the computational cost from blowing up when the number of discretization cells increases, the weights, $t_k$, should be nonzero only locally around the face. In the most extreme case, one can choose to approximate the gradient using the pressure values from the two cells immediately next to the face to compute the discharge from $c_i$ to $c_j$:

$$u_{ij} := \int_{S_{ij}} \mathbf{u} \cdot \mathbf{n} dS = -\int_{S_{ij}} \nabla p \cdot \mathbf{K} \cdot \mathbf{n} dS \approx -t_{ij}(p_j - p_i). \qquad (8)$$

In this case, the face transmissibility is computed as the harmonic average of the half transmissibilities:

$$t_{ij} = \frac{\alpha_i \alpha_j}{\alpha_i + \alpha_j}, \qquad (9)$$

which in turn are given as

$$\alpha_i = \frac{A_{ij} \mathbf{n}_{ij} \cdot \mathbf{K}_i}{\mathbf{d}_{ij} \cdot \mathbf{d}_{ij}} \cdot \mathbf{d}_{ij}. \qquad (10)$$

Here, the aperture-weighted face area is denoted by $A_{ij}$ and the unit normal vector $\mathbf{n}_{ij}$ on the face points outward from $c_i$. The permeability tensor and the distance vector between the cell centre $x_{cell}$ and the face centre $x_f$, $\mathbf{d}_{ij} = x_f - x_{cell}$, also belong to $c_i$. In the left part of Figure 2, the quantities needed for calculating the half transmissibility $\alpha_i$ of the left cell and centre face are shown. The TPFA is described in more detail in (Aziz & Settari, 1979) and is the standard discretization in commercial reservoir simulation.

### 2.3.2 Multi-Point Flux Approximation

While straightforward to derive and implement and computationally very efficient due to matrix sparsity, the TPFA scheme is only consistent when the grid is aligned with the principal directions of the permeability tensor. Defining the normal component of the permeability as $\mathbf{w}_{ij} = \mathbf{K}_i \cdot \mathbf{n}_{ij}$, the requirement is that $\mathbf{w}_{ij}$ and $\mathbf{d}_{ij}$ be aligned for all cell-face combinations. If this is not the case, as in the right part of Figure 2, a more sophisticated flux approximation is called for.

To construct the multi-point flux approximation, we split each face at the face centre. The transmissibility calculations are done for all the global sub-faces, and then summed for all faces and, finally, for all cells according to Equation (5).

The sub-face transmissibility computations are performed on local systems defined by an O-shaped interaction region around each grid node. The interaction region connects the surrounding cell and face centroids. Pressure continuity is enforced at a single point on each sub-face, given as $x_c = (1 - \eta)x_f + \eta x_n$, with $x_c$ denoting the continuity point and $x_n$ the vertex associated with the sub-face; see Figure 2. We set $\eta = 0$ on Cartesian grids, and $\eta = 1/3$ on simplex grids, as suggested in (Aavatsmark, 2002; Edwards & Rogers, 1998; Klausen, Radu, & Eigestad, 2008; Friis, Edwards, & Mykkeltveit, 2009), where more details on the MPFA method are to be found.

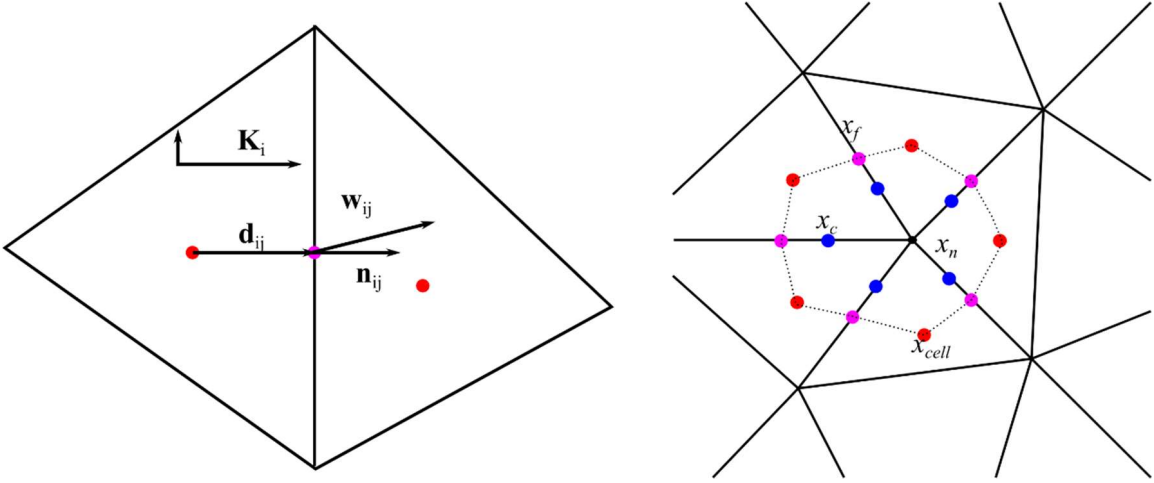

*Figure 2. To the left, entities for TPFA transmissibility computations between the left cell i and the right cell j. To the right, the O-shaped interaction region used for MPFA transmissibility calculations. Cell centres in red, face centres in purple, continuity points in blue and the central vertex of the interaction region in black.*

Assuming the pressure to be linear within each cell of the region facilitates the gradient evaluation needed to adhere to Darcy's law. The flux from $c_i$ over sub-face $j$ of area $A_j$ becomes

$$u_j^i = \sum_k A_j \mathbf{n}_j \cdot \mathbf{K}_i p_k \nabla \phi_k(x), \tag{11}$$

where $\phi_k$ are the linear basis functions of the interaction region and the summation index $k$ includes the cell centre of $c_i$ and all continuity points. One such relation is obtained for each face-cell combination of the interaction region. By equating them over each sub-face (continuity of flux), the intermediate continuity point pressures may be eliminated. All sub-face fluxes of the interaction region are assembled for each interaction region yielding

$$\mathbf{u} = \mathbf{T}\mathbf{p}. \tag{12}$$

Here, $\mathbf{p}$ consists of the cell pressures of the interaction region and $\mathbf{T}$ is the matrix of sub-face transmissibilities accounting for the local permeability and geometry of the interaction region. These transmissibilities are assembled to the global matrix of Equation (7).

### 2.3.3 Transport discretization

The transport equation is also integrated over each cell to yield a FV formulation. The upwind term is integrated by parts and then discretized with upwind concentration evaluation, i.e., for the advection between $c_i$ and $c_j$ we use

$$T_{upw} = \begin{cases} T_i & \text{if } u_{ij} > 0, \\ T_j & \text{if } u_{ij} < 0. \end{cases} \tag{13}$$

Assuming stationary fluxes, sinks and sources, implicit Euler temporal discretization is applied, yielding for time step $n$

$$V_i \frac{T_i^n - T_i^{n-1}}{\Delta t} - \sum_j A_{ij} T_{upw}^n u_{ij} = s_i \qquad (14)$$

for the concentration $T_i$ at the centre of $c_i$ of volume $V_i$.

## 2.4 Subdomain Coupling

We first consider the coupling of the flux discretizations between different subdomains. In the discretizations of the individual subdomains, we impose Neumann conditions for all boundaries internal to the full domain. This allows us to couple the subdomains by simply discretizing the flux over all faces of each internal boundary.

We consider the coupling of two subdomains $\Omega_i^{d+1}$ and $\Omega_j^d$ over the common boundary $\Gamma_{ij}$. Letting the subscript indexes $i$ and $j$ also indicate cells from the respective subdomains, we discretize the flux from $c_i \in \Omega_i^{d+1}$ to $c_j \in \Omega_j^d$ according to Karimi-Fard et al. (2004):

$$\hat{\mathbf{n}}_{ij} \cdot \hat{\mathbf{u}}_{ij} = -t_{ij}(p_j - \hat{p}_i), \qquad (15)$$

with $t_{ij} = \frac{\hat{\alpha}_i \alpha_j}{\hat{\alpha}_i + \alpha_j}$. For the half face transmissibility calculation, the lower dimensional subdomain is extended by half an aperture $a$ in the normal direction; hence

$$\alpha_j = \frac{\mathbf{n}_{ji} \cdot \mathbf{K}_j}{a/2} A_{ij} \qquad (16)$$

on the fracture side of the interface. On the matrix side, $\hat{\alpha}_i$ is computed according to Equation (10). The DFM approach inherently accounts doubly for the fracture domains both in cell volumes and distance computations, and so relies on $a \ll h_{min}$. The distance inconsistency can be corrected by using the adjusted distance vector

$$\mathbf{d}_{ij}^* = \mathbf{d}_{ij}\left(1 - \frac{a}{2|\mathbf{d}_{ij}|}\right); \qquad (17)$$

see e.g. Sandve, Nodbotten and Berre (2012).

As the flux over the face may be expressed in terms of one pressure value from each of the two subdomains only, the couplings remain entirely local and can be computed independently. The coupling conditions are arranged in four blocks and assembled into the global solution matrix, as illustrated in Figure 1. Thus, the diagonal blocks consist of the sum of the internal discretization and the contribution from the subdomain coupling. The off-diagonal blocks will be nonzero only where the corresponding subdomains interact.

For the transport problem, upwind coupling terms are calculated just as in the internal discretizations. For the area, the face area of the higher-dimensional cell is used.

## 2.5 Intersection cell removal

In realistic fractured porous domains, the fracture aperture is commonly several orders of magnitude smaller than the domain length scale. This implies high contrasts in discretization cell size, as the smallest cell volumes in dimension $d$ typically scale according to

$$V(c) \propto h_{min}^d \, a^{N-d}. \tag{18}$$

Strong cell size variation may affect the solution quality and efficiency of an iterative linear solver through high condition number $C = cond(\mathbf{A})$ of the matrix $\mathbf{A}$ from Equation 7. As pointed out and demonstrated by Karimi-Fard, Durlofsky and Aziz and Walton, Unger, Ioannidis and Parker (2004; 2017), removing the one- and zero-dimensional cells at fracture intersections may improve the condition number significantly, while also relaxing time step restrictions for transport and multiphase flow simulations.

Karimi-Fard et al. (2004) proposed an elimination procedure for the TPFA, known as a Star-Delta transformation, where direct connections are introduced between the higher-dimensional cells. These connections are obtained by adding $n_f$ fluxes between all the $n$ higher-dimensional cells meeting at the particular intersection, see Figure 3, where $n = 4$ leads to six cell-pair combinations and hence $n_f = 6$ new fluxes. The transmissibility for each cell pair at the intersection is computed as

$$T_{ij} = \frac{\alpha_i \alpha_j}{\sum_k^n \alpha_k}, \tag{19}$$

where $\alpha_k$ is computed as in Equation 10. This popular transformation technique was extended to the MPFA method by Sandve et al. (2012).

By ignoring the intersection cell, the Star-Delta elimination in effect assigns infinite permeability to the intersection. When the permeability of the intersection is in fact low, which could be the case if a blocking fracture crosses a permeable one, this may lead to significant overestimation of fracture connectivity, as demonstrated in Section 3.1 and further studied in (Stefansson, 2016).

To remedy this shortcoming, we here introduce a new technique for intersection cell removal based on the Schur complement; see, for example, (Zhang, 2005). We first discretize the problem including intersection cells and then identify the degrees of freedom to be <u>k</u>ept and those to be <u>e</u>liminated. Associating the matrix rows and columns corresponding to the two sets with subscripts $k$ and $e$, respectively, we obtain the reduced system matrix

$$\mathbf{A}_\mathrm{r} = \mathbf{A}_\mathrm{kk} - \mathbf{A}_\mathrm{ke} \mathbf{A}_\mathrm{ee}^{-1} \mathbf{A}_\mathrm{ek}. \tag{20}$$

The reduced right hand side is computed similarly, yielding the reduced system $\mathbf{A}_\mathrm{r} \mathbf{p}_\mathrm{r} = \mathbf{b}_\mathrm{r}$. Solving this system for the kept pressures, we accurately account for intersections of fractures of arbitrary permeability ratios while avoiding the minute volumes of the removed cells. The technique yields accurate pressure solutions also in cases involving flow along the intersection when 1D cells are eliminated. Further, the Schur complement reduces to the above mentioned Star-Delta transformation in the limiting cases of infinite normal permeability and zero tangential permeability of the eliminated cells, as demonstrated in Section 3.1.1. Note that this elimination technique is not particular to FV discretizations.

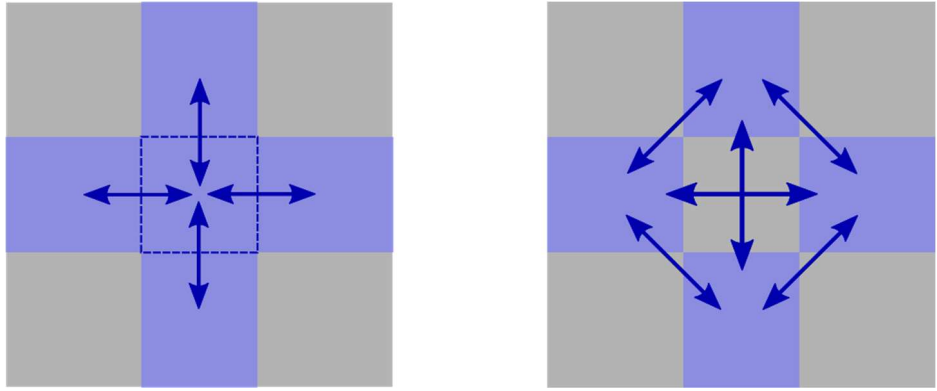

*Figure 3. Conceptual sketch of the fluxes at a fracture intersection with and without an intersection cell, left and right, respectively.*

If fluxes between the grid cells are to be computed, there are two ways to proceed after the elimination. The first option is to back-compute the pressure in the intersection cells and proceed as if the cells were never removed. While straightforward using the Schur complement, this is in many situations somewhat counterintuitive. If transport simulations are performed using the resulting flux field, the original time step restrictions have to be honoured.

Not reintroducing the lower-dimensional cells calls for flux computations between the higher-dimensional neighbours directly, see Figure 3. This may be done directly from the reduced discretization, in the same manner as described in Section 2.3 for the internal subdomain fluxes. This option is not entirely unproblematic, especially if 1D cells are eliminated. Flow along the intersection will be distributed more evenly between the intersecting fractures compared to a non-eliminated or back-calculated solution, leading to additional "leakage" from the main flow path to less permeable regions, as shown by Sandve, Nordbotten and Berre (2012).

## 2.6 Anisotropy

As mentioned above, anisotropic permeability calls for more sophisticated discretizations such as the MPFA. Using a DFM model, the background matrix permeability may also include the upscaled contribution from fine-scale fractures, which may enhance the matrix anisotropy.

Using TPFA for the fracture-matrix coupling might impair the solution quality in the case of matrix anisotropy as the half-face transmissibility on the matrix side, $\hat{\alpha}_i$, will be computed inconsistently, cf. Section 2.4. This is investigated in Section 3.2.

As shown by Mokhtari et al. (2015) and Tan et al. (2017), fracture planes may also exhibit intrinsic permeability anisotropy; that is, the permeability depends on the direction considered in the fracture plane itself. Given that the fracture network and its connectivity oftentimes dominate the flow through the domain, it is conceivable that fracture anisotropy may substantially influence the overall behaviour of the system. We investigate the impact of fracture anisotropy in Section 3.3.

# 3 Results

The results are presented in the form of plots of pressure and tracer distributions and plots of locally monitored tracer concentrations throughout the transport simulation time interval. As the main purpose of the plots is to qualitatively illustrate the solutions and differences among them, we omit colour bars, noting that scales are consistent and all colour maps use blue for low values and red for high ones. We report errors of both pressure and final tracer fields as compared to various reference solutions. The errors are computed either over the entire domain or over specified subdomains, always in the discrete $L^2$ norm.

The numbering of the test cases corresponds to the subsections, so that Case 1.1 appears in Section 3.1.1 etc. We reiterate that source code is available at www.github.com/pmgbergen/porepy, as is additional visualization material intended as a supplement to the test case descriptions, figures and tables provided in this section.

## 3.1 Intersection elimination

This section is devoted to an examination of the Schur complement procedure for elimination of fracture intersections introduced in Section 2.5. We compare the pressure solution to a solution based on a reference discretization where the intersection cells have been kept and a solution where the Star-Delta transformation is used for elimination of intersection cells. In addition to assessing the solution quality through pressure and tracer field comparisons, we report the condition numbers of the respective solution matrices. An indication on the improvement offered by the eliminations is provided by the ratio between the condition number of the full solution and the eliminated, which we denote by $R_C$. For these test cases all permeabilities are isotropic, and we thus only apply TPFA.

### 3.1.1 Elimination of a 0D intersection

We start by a systematic investigation of the elimination procedures in a synthetic one-dimensional test case containing a single intersection. We incrementally change the fracture permeabilities to reveal when the Schur complement and Star-Delta procedures differ.

The test case geometry is shown in Figure 4. The domain contains two orthogonal fractures, one normal to and one parallel to the pressure gradient, both having an aperture $a = 10^{-2}$. The matrix permeability is set to unity, whereas the permeabilities of the fractures vary between $10^{-3}$ and $10^3$ to yield a range of different ratios between the permeability of the horizontal (gradient-aligned) fracture, $K = K_h I$, and that of the vertical fracture, $K = K_v I$. The permeability of the 0D intersection cell is inherited from the vertical fracture. Homogeneous Neumann conditions are imposed on the horizontal boundaries, and the Dirichlet condition $p_D = 1 - x$ is imposed on the vertical boundaries.

In Figure 4, the pressure plots for the permeability ratio yielding the largest difference between the elimination techniques are shown. The pressure jump across the vertical fracture is markedly higher for the reference and Schur complement solution than for the Star-Delta one, where the blocking effect at the intersection is lost. We include the pressure plot of a simulation where the intersection cell permeability is $10^{10}$. The solution is practically identical to the results using the Star-Delta transformation, which illustrates how the latter is the limiting case of high intersection permeability, as discussed in Section 2.5.

A detailed report of the pressure errors and condition numbers is given in Figure 5. We see machine precision agreement between the solution for which the Schur complement technique was applied and the non-reduced solution. The solution where the Star-Delta transformation was used, however, breaks down when the $\nabla p$-aligned fracture is permeable and the other one is not. Under these conditions, the horizontal fracture strongly influences the solution, but the blocking effect of the vertical fracture is not captured by the Star-Delta transformation. The two fracture elimination schemes bring about significant and very similar condition number improvements, up to more than a factor of $10^2$. Note that while the magnitude of $R_c$ in general depends on the aperture and permeability values in question, Figure 5 indicates that the improvement is larger when low-permeable intersections are eliminated.

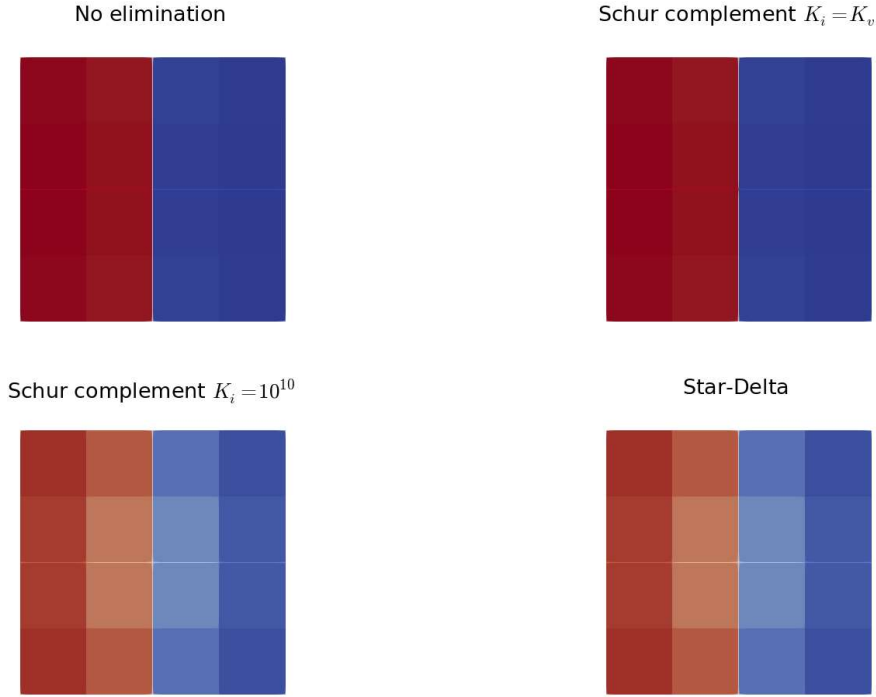

*Figure 4. Pressure fields of Case 1.1, $K_h = 10^3$ and $K_v = 10^{-3}$. The case is designed so that the intersection permeability, $K_i$, equals that of the vertical fracture. The bottom left Schur complement solution with $K_i = 10^{10}$ is included for comparison to the Star-Delta solution.*

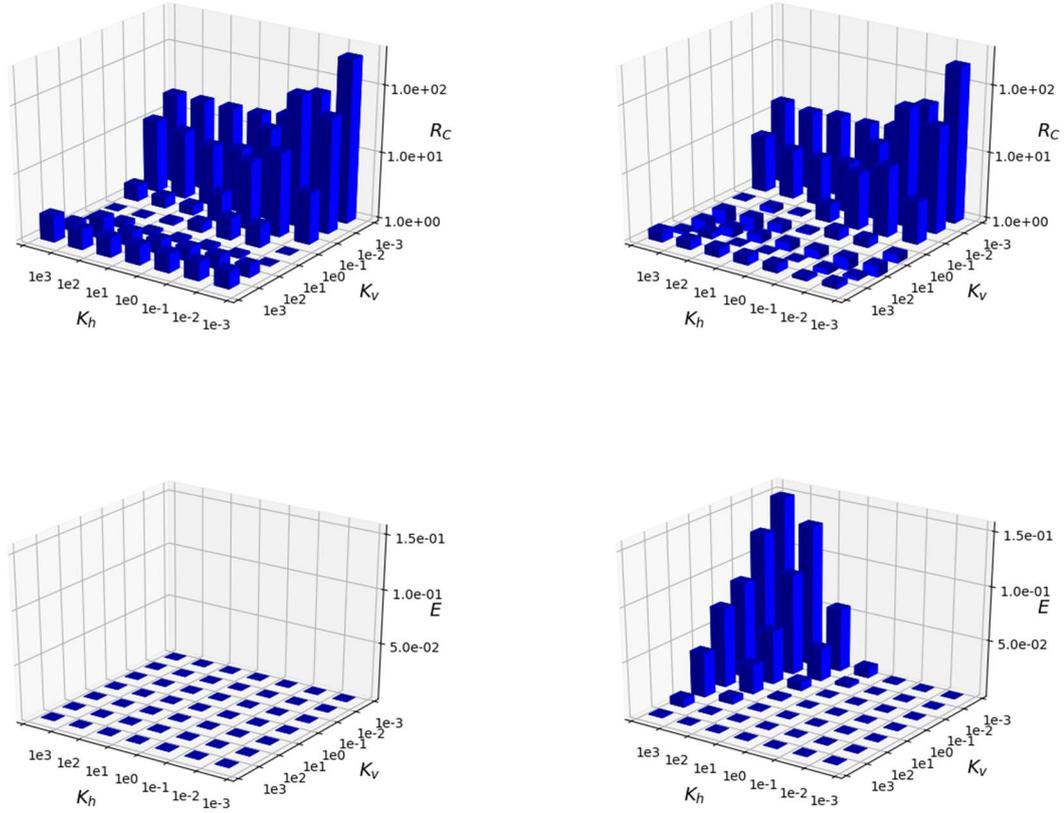

*Figure 5. Condition number ratios $R_C$ (top) and pressure errors (bottom) for Case 1.1 for different values of the permeabilities of the horizontal and vertical fracture. Schur complement and Star-Delta elimination were used for the sub-figures to the left and to the right, respectively.*

### 3.1.2 Elimination of multiple 0D intersections

We now consider a 2D test case presented in the benchmark study (Flemisch, et al., 2018). In one of the test cases of that study containing blocking intersections in conducting fractures, applying the Star-Delta intersection cell elimination procedure to the TPFA and MPFA considered in the present work was shown to decrease solution quality. We here investigate the Schur complement elimination procedure for that particular test case.

The test case contains ten fractures of two different permeabilities. The intersection permeabilities are set to the harmonic average of the crossing fractures, making them in effect impermeable. For further details on the set-up we refer to (Flemisch, et al., 2018), where this benchmark test case appears in section 4.3.2. The grid used here for the simulation with the Schur complement is slightly different from the one used for the other three. However, the number of fracture and matrix cells are almost identical.

We display two of the solutions in Figure 6 and the errors and characteristics of the matrices in Table 1. The plots indicate that using the Star-Delta technique, significant flow passes through the blocking fractures where they intersect with conducting ones. The errors of the solutions obtained using the Schur complement elimination are very similar to those obtained without elimination and considerably lower than those where the Star-Delta elimination were, and thus indicate that the Schur complement technique does not suffer from the restrictions of the Star-Delta transformation related to blocking intersections. Moreover, the condition number of the Schur complement matrices are slightly lower than for the two Star-Delta matrices. This means that we obtain the full benefits of the Star-Delta technique while retaining solution quality.

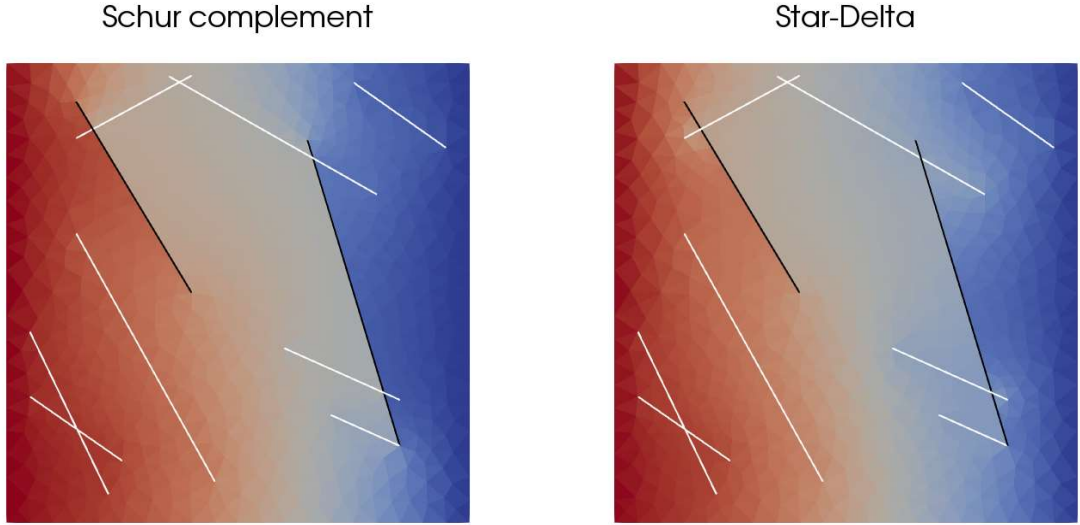

*Figure 6. The pressure distributions for Case 1.2 obtained using the two elimination techniques with MPFA discretization. Permeable and impermeable fractures are indicated in white and black, respectively.*

*Table 1. Errors and matrix characteristics for Case 1.2.*

| Method | Fracture error | Matrix error | Matrix density | Condition number |
|---|---|---|---|---|
| TPFA, no elimination | $1.3 \cdot 10^{-2}$ | $1.3 \cdot 10^{-2}$ | $2.8 \cdot 10^{-3}$ | $2.0 \cdot 10^5$ |
| TPFA, Schur complement | $1.4 \cdot 10^{-2}$ | $1.4 \cdot 10^{-2}$ | $2.8 \cdot 10^{-3}$ | $2.1 \cdot 10^4$ |
| TPFA, Star-Delta | $5.1 \cdot 10^{-2}$ | $6.7 \cdot 10^{-2}$ | $2.8 \cdot 10^{-3}$ | $3.1 \cdot 10^4$ |
| MPFA, Schur complement | $1.2 \cdot 10^{-2}$ | $1.4 \cdot 10^{-2}$ | $8.0 \cdot 10^{-3}$ | $2.0 \cdot 10^4$ |
| MPFA, Star-Delta | $5.1 \cdot 10^{-2}$ | $6.7 \cdot 10^{-2}$ | $8.5 \cdot 10^{-3}$ | $3.1 \cdot 10^4$ |

### 3.1.3   Elimination of a 1D intersection

We now turn to a 3D case where the 1D intersection is eliminated. In general, this is inherently different from elimination of 0D intersections, as the tangential flow along the intersection is affected. In particular, flux calculations become more involved, as discussed in Section 2.5.

The unit cube domain of unitary permeability contains two fractures lying in the *xy*- and *yz*-plane, respectively. Both have aperture $10^{-6}$, and the isotropic fracture permeabilities are $10^6$ and $10^{-6}$, respectively. The boundary conditions are homogeneous Neumann, except for the pressure conditions $p_D(x = 0) = 1$ and $p_D(x = 1) = 0$. Similarly, tracer concentration values of zero are prescribed at all external boundaries, whereas the initial tracer concentration is one throughout the domain. The transport simulation runs to $t = 0.5$.

We compare three solutions with TPFA in the entire domain: one where no elimination of intersection cells is performed, one where the Schur complement elimination is used and one using the Star-Delta elimination. The errors of the eliminated solutions compared to the full one are listed in Table 2. The tracer concentrations shown in Figure 7 show qualitative differences: If the blocking nature of the intersection is honoured, the tracer is largely forced around the vertical fracture, whereas it shoots through in the case of the Star-Delta elimination. The perfect agreement between the solution without elimination and the Schur complement, and difference to the Star-Delta, are highlighted by the concentration time series in Figure 8, recorded at the centre of the out-flow boundary.

For the solutions where the 1D intersection cells were eliminated, the fluxes at the intersections are computed directly between the neighbouring 2D cells, as described in Section 2.5 and indicated to the right in Figure 3. The very low error of the Schur complement solution illustrates that this procedure is adequate also for flux computations provided that the flow along the fracture is negligible.

*Table 2. Errors of the eliminated solutions compared to the full solution (with 1D cells) for Case 1.3*

| Elimination | Pressure error | Tracer error | Condition number improvement |
|---|---|---|---|
| Schur complement | $2.2 \cdot 10^{-15}$ | $1.5 \cdot 10^{-12}$ | $6.0 \cdot 10^3$ |
| Star-Delta | $3.2 \cdot 10^{-2}$ | $1.2 \cdot 10^{-1}$ | $4.1 \cdot 10^3$ |

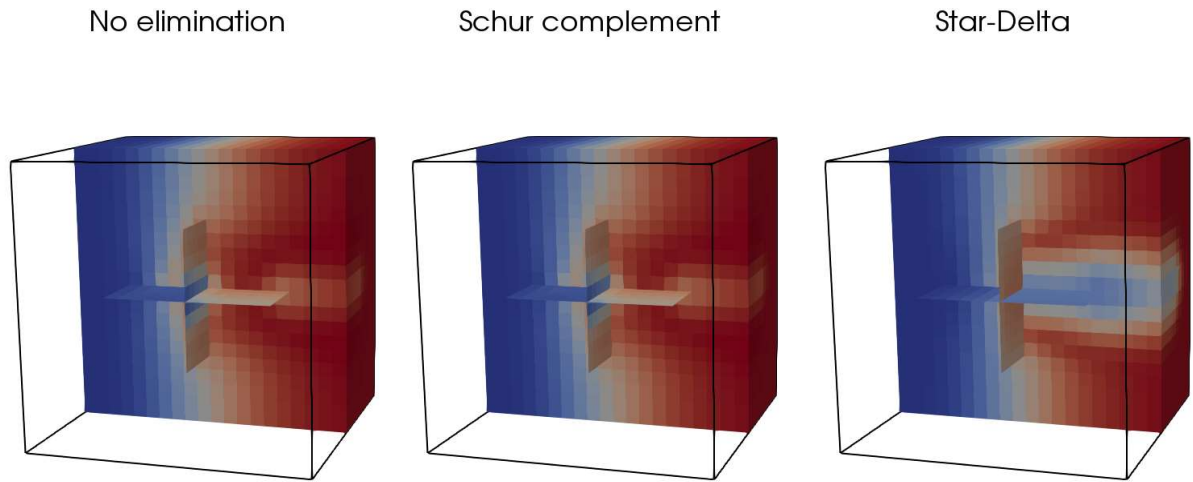

*Figure 7: Case 1.3 tracer concentrations at time $t = 0.5$.*

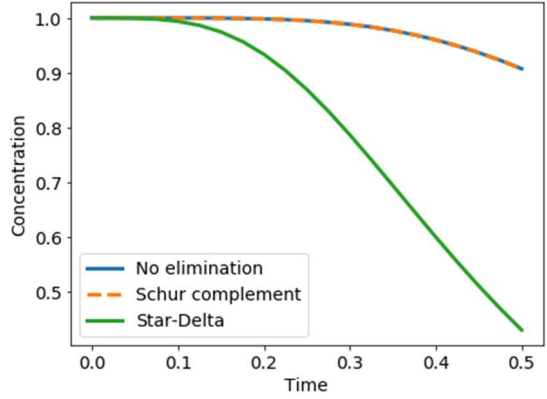

*Figure 8. Monitored tracer concentration at the centre of the outflow boundary throughout the transport simulation for Case 1.3.*

## 3.2 Subdomain coupling with anisotropy in matrix permeability

Turning back to the unit square domain, we investigate the TPFA coupling's behaviour in face of anisotropy in the higher dimension by considering a single-fracture case. In particular, we investigate the impact of non-alignment of the distance vector $\hat{\mathbf{d}}_{ij}$ and normal permeability component $\hat{\mathbf{w}}_{ij}$ on the matrix side of the subdomain interface, with subscript indexes as in Section 2.4. This means that we apply the TPFA to a problem for which it is known to be inconsistent, but only in a small part of the domain.

Setting the matrix permeability to unity, we include one highly permeable ($\mathbf{K} = 10^4\,\mathbf{I}$) horizontal fracture of aperture $10^{-3}$ at $y = 0.5$. No-flow conditions are set for both pressure and tracer except at two diametrically opposing corners of the domain, where Dirichlet conditions are enforced. A series of simulations are performed with varying anisotropy ratio $K_{max}/K_{min}$ and angle (10°, 30° and 60° relative to the coordinate axes) in the matrix, i.e.,

$$\mathbf{K} = \mathbf{R}(\theta) \begin{bmatrix} K_{max} & 0 \\ 0 & K_{min} \end{bmatrix},$$

where $\mathbf{R}$ is the rotational matrix and $\theta$ the angle of rotation relative to the $x$ axis.

We compare the reduced model using the TPFA coupling and MPFA for internal discretizations to a equi-dimensional model using MPFA in the entire domain in terms of both pressure and tracer errors. The reference solution is equi-dimensional and computed on a fine mesh ($256 \times 256$ cells). We refine the coarse solutions from $4 \times 4$ to $32 \times 32$ cells. Because the matrix cell size $h_{min}$ is only one order of magnitude larger than the aperture for the finest grid, we use the aperture corrected distance vector introduced in Section 2.4, thus avoiding that the inconsistency of the hybrid-dimensional approach obscures the convergence comparison.

In Figure 9, we show the results for the 30° rotation (the two other angles yielded no qualitative differences). The results indicate that the TPFA coupling performs very acceptably. There is, however, one significant limitation to the conclusion. For the anisotropy, we consider a permeability ratio range of $1 \leq K_{max}/K_{min} \leq 6$. For larger values the results showed oscillations in the matrix pressure away from the fracture. We attribute this to monotonicity issues for the MPFA unrelated to the presence of fractures, see (Nordbotten, Aavatsmark, & Eigstad, 2007). We expect that the errors related to the TPFA coupling scheme are enhanced for larger anisotropy ratios, although explicitly accounting for the large fractures relieves us from the most severe anisotropies characteristic of effective permeabilities in the context of fracture upscaling.

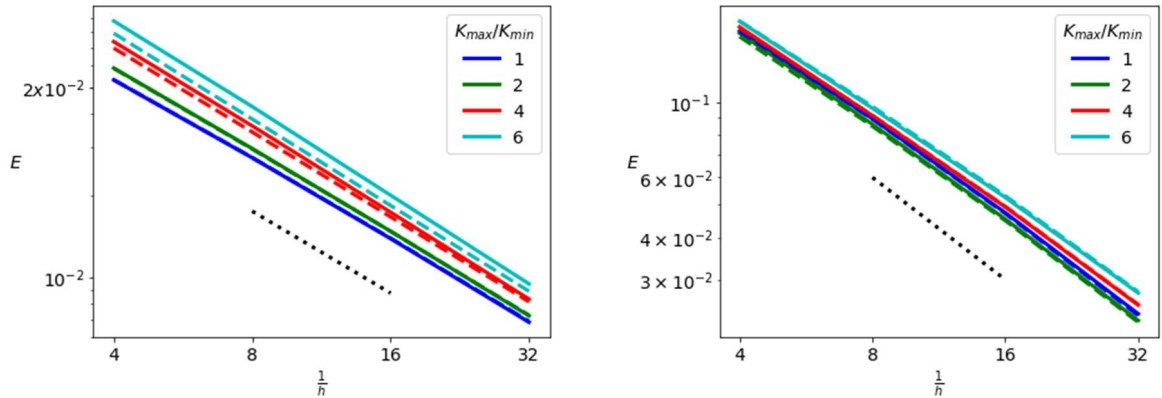

Figure 9. Pressure and tracer error plots for Case 2, left and right, respectively. Line colours represent different permeability ratios. Solid lines correspond to the reduced model, dashed lines to the equi-dimensional MPFA. Dotted black lines indicate a linear slope.

### 3.3 Anisotropy in fracture permeability

A 3D test case investigating the effect of permeability anisotropy in a 2D fracture is now studied. We explore whether the anisotropy calls for the use of MPFA or if TPFA is admissible. We again consider a case with limited geometric complexity and assign isotropic permeability in the matrix. To further isolate the effect of fracture anisotropy, we use a Cartesian grid for our unit cube domain, ensuring alignment of face normals and distance vectors in the matrix, cf. Section 2.3.1.

The domain contains a single, horizontal fracture of aperture $10^{-3}$ cutting the domain at $z = 0.5$. No-flow conditions are imposed on all but a part of the top and bottom boundaries, for which we set $p_D = 0$ and $p_D = 1$, respectively. With unit permeability in the matrix and a fracture permeability tensor with principal permeabilities of $10^3$, $1/3 \cdot 10^3$ and $10^3$ rotated 45° in the fracture plane:

$$\mathbf{K} = \mathbf{R}(45°) \begin{bmatrix} 10^3 & 0 & 0 \\ 0 & 1/3 \cdot 10^3 & 0 \\ 0 & 0 & 10^3 \end{bmatrix} = \begin{bmatrix} 2/3 \cdot 10^3 & -1/3 \cdot 10^3 & 0 \\ -1/3 \cdot 10^3 & 2/3 \cdot 10^3 & 0 \\ 0 & 0 & 10^3 \end{bmatrix},$$

we obtain a system strongly influenced by the fracture flow. Again, the initial and boundary condition for the tracer is $T(t = 0) = 1$ in $\Omega$ and $T = 0$ on $\delta\Omega$. The transport simulation runs to $t = 30$.

Comparison of the pressure solutions obtained using TPFA and MPFA yields the relative errors reported in Table 3, revealing only slight differences. The tracer solutions are compared both in terms of $L^2$ errors and a time series of the concentration at the outflow boundary in Figure 10. By contrast, these reveal disparities between the flow field of significant magnitude.

As mentioned in Section 1, we may also combine the discretizations on the different subdomains. Exploiting that the permeability is anisotropic only in the fracture, we construct a hybrid discretization with TPFA in the matrix and MPFA in the fracture, coupled together with the TPFA coupling. As evidenced from the tracer results, this discretization is sufficient to capture the anisotropy. We also report the CPU times obtained with our implementation for the discretization of the problem in Table 3. Although these should as always be treated with caution, they indicate that the computational cost of the hybrid discretization is substantially reduced compared to the MPFA, although somewhat larger than for the TPFA.

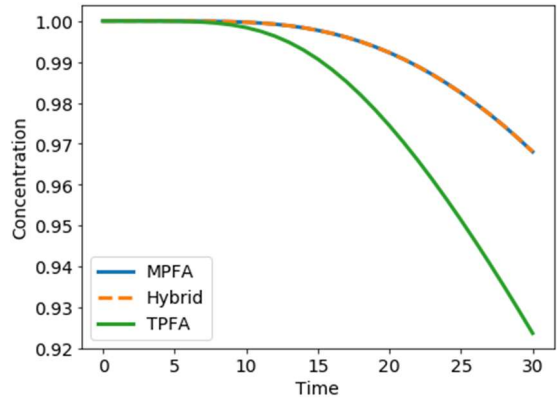

*Figure 10. Time series of tracer monitored in the cell of the outflow boundary, Case 3.*

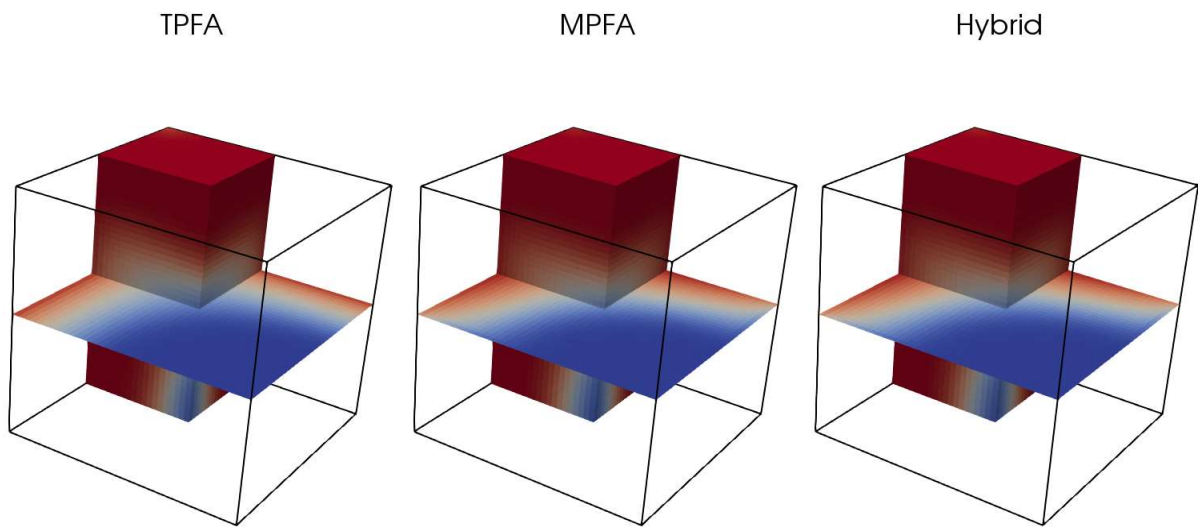

*Figure 11. Tracer distributions at $t = 30$ for Case 3. There is injection in the top, back corner and production takes place at the opposite corner. The 2D fracture plane is left uncut for visualization purposes.*

*Table 3. Difference to the MPFA reference in solutions and discretization times for Case 3.*

| Flux discretization | Subdomain | Difference to MPFA pressure | Difference to MPFA tracer | Discretization time, relative to MPFA |
|---|---|---|---|---|
| TPFA | 3D matrix | $5.0 \cdot 10^{-4}$ | $5.8 \cdot 10^{-2}$ | $3.3 \cdot 10^{-2}$ |
|  | 2D fracture | $8.4 \cdot 10^{-4}$ | $1.6 \cdot 10^{-1}$ |  |
| Hybrid | 3D matrix | $1.8 \cdot 10^{-10}$ | $8.6 \cdot 10^{-8}$ | $8.3 \cdot 10^{-2}$ |
|  | 2D fracture | $2.7 \cdot 10^{-10}$ | $2.4 \cdot 10^{-7}$ |  |

## 3.4 Matrix heterogeneity

In the final test case, we provide an example of how the hierarchical subdomain approach and coupling flexibility may be advantageously applied. We investigate a 3D case displaying both strong fracture influence and matrix flow with a fracture geometry containing several different intersection configurations.

We choose a conductive fracture network, with $\mathbf{K} = 10^5 \mathbf{I}$ and $a = 10^{-6}$, symmetric about the $xy$-plane and blocked by one fracture, with $\mathbf{K} = 10^{-5}\mathbf{I}$, particularly the circular one at the centre of the domain as shown in Figure 12. Intersection permeabilities are inherited from the least permeable of the intersecting fractures. The network is surrounded by a matrix of permeability $\mathbf{K} = 10^{-2}\mathbf{I}$ in the upper half and $\mathbf{K} = 10^{-3}\mathbf{I}$ in the lower half. No-flow conditions are imposed on the vertical boundaries and $p, T = 0$ at the horizontal ones, whereas a well injection with a tracer concentration of one is located in a fracture intersection cell near the centre of the domain. Starting from an initial tracer solution of $T = 0$, the simulation runs until $t = 2$.

Since all permeabilities are isotropic, we discretize using TPFA in all dimensions. In this test case, we apply the Schur complement intersection elimination procedure, and compare the solution with one where the intersection cells were kept. The solution is visualized in Figure 12. The tracer distribution clearly shows how the symmetry of the test case is broken by the heterogeneous matrix permeability. In the upper half, a substantial part of the flow takes place in the matrix, whereas virtually all of it is conducted through the fractures in the lower half. Again, considering the solution without the elimination as reference, negligible errors and considerable condition number improvement are observed for the Schur complement procedure, see Table 4.

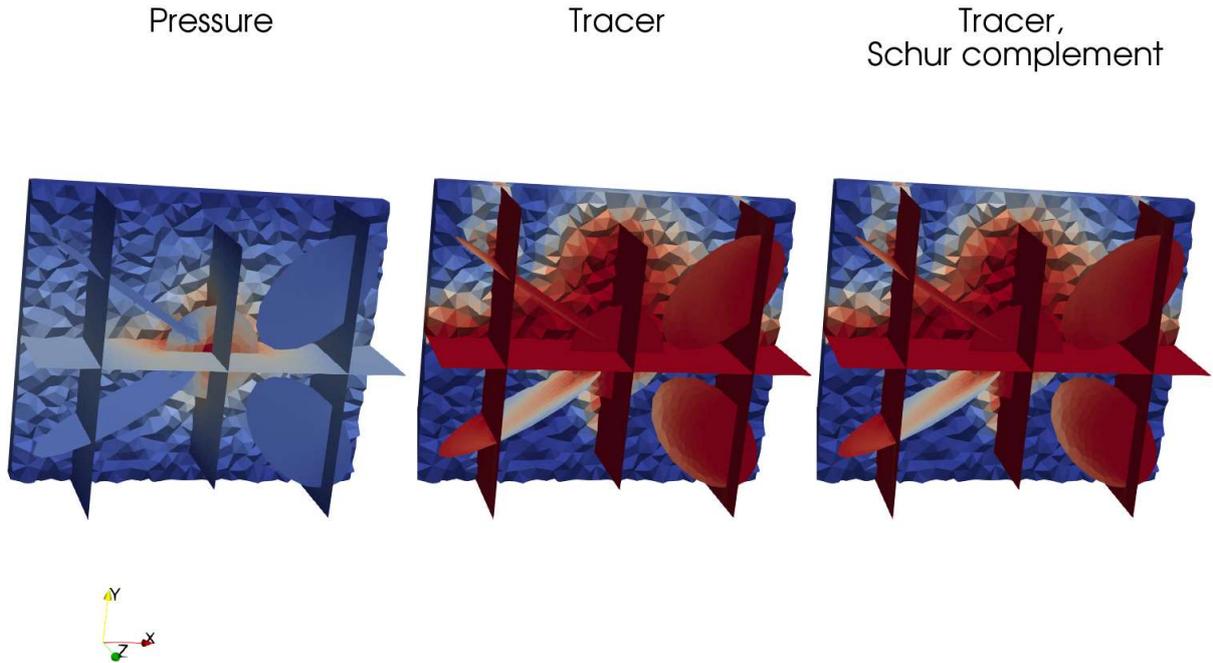

*Figure 12: Distribution of pressure (left) and tracer at $t = 2$ for a solution obtained without intersection cell removal (middle) and one where the Schur complement procedure was applied (right) for Case 4. The injection well placement may be inferred from the pressure plot.*

*Table 4. Errors and condition number improvement for the elimination procedure applied in Case 4.*

| Elimination | Pressure error | Tracer error | Condition number improvement |
|---|---|---|---|
| Schur complement | $1.6 \cdot 10^{-11}$ | $4.8 \cdot 10^{-11}$ | $8.6 \cdot 10^3$ |

# 4 Conclusion

Choosing a hybrid-dimensional conforming DFM model for simulations of flow in fractured porous media may relieve some of the challenges posed by the fractures. We have here considered procedures to alleviate two of the remaining restrictions, namely the coupling of fractures and rock matrix and the wide range in discretization cell size. Firstly, we demonstrated how discretizing the subdomains individually and coupling them two at a time using a hierarchical framework for mixed-dimensional media problems simplifies the implementation and offers considerable modelling and discretization flexibility, also facilitating reuse of code. Secondly, we proposed a Schur complement technique to elimination of intersection cells with a considerably broader scope than existing techniques. The technique is not limited to FV discretizations, and may be applied to both 0D and 1D intersections of arbitrary permeabilities without incurring additional pressure error, while the computation of fluxes introduces some issues in cases involving tangential flow along 1D intersections.

The numerical test cases presented demonstrate the improved accuracy of the proposed intersection cell elimination technique compared to existing solutions. Further, they show that the proposed TPFA coupling between fractures and matrix does not lead to severely reduced solution quality even in the presence of moderate anisotropy. We also demonstrate how parameter heterogeneities and anisotropies have implications on the appropriateness of different models and discretizations, and show how the inherent flexibility of our approach renders it suitable for simulations considering heterogeneous parameter regimes. The approach retains the benefits of FV discretizations and renders it feasible for simulation of various processes in fractured porous media.

# 5  Acknowledgements

The work has partly been funded by the Research Council of Norway through grant number 267908.